\documentclass[twoside,11pt,reqno]{amsart}
\usepackage{amsfonts, amsbsy, amsmath, amsthm, amssymb, latexsym, verbatim, enumerate}

\usepackage[english]{babel}

\newtheorem{propo}{Proposition}[section]

\newtheorem{lemma}[propo]{Lemma}
\newtheorem{corol}[propo]{Corollary}

\newtheorem{theo}[propo]{Theorem}

\theoremstyle{definition}

\newcommand{\ld}{,\ldots ,}

\newcommand{\ra}{ \rightarrow }

\newcommand{\lan}{ \langle }
\newcommand{\ran}{ \rangle }

\newcommand{\diag}{\mathop{\rm diag}\nolimits}

\newcommand{\Id}{\mathop{\rm Id}\nolimits}

\newcommand{\Om}{\Omega}

\newcommand{\FF}{\mathop{\Bbb F}\nolimits}

\newcommand{\QQ}{\mathop{\Bbb Q}\nolimits}

\newcommand{\ZZ}{\mathop{\Bbb Z}\nolimits}

\newcommand{\al}{\alpha}

\newcommand{\ep}{\varepsilon}

\newcommand{\lam}{\lambda }

\newcommand{\up}{^{-1}}

\newcommand{\si}{\sigma }

\def\11{$(1)$}
\def\22{$(2)$}
\def\33{$(3)$}

\def\d12{{_{12}}}

\def\acf{{algebraically closed field }}

\def\au{{automorphism }}

\def\ei{{eigenvalue }}
\def\eis{{eigenvalues }}

\def\f{{following }}

\def\ho{{homomorphism }}

\def\ii{{if and only if }}

\def\ir{{irreducible }}

\def\irt{{irreducible. }}
\def\irr{{irreducible representation }}

\def\itf{{It follows that }}

\def\mult{{multiplicity }}

\def\rep{{representation }}
\def\reps{{representations }}

\def\F{{\rm F}}

\newcommand{\om}{\omega}
\newcommand{\el}{\end{lemma}}

\def\ag{algebraic group }

\newcommand{\bl}{\begin{lemma}\label}
\newcommand{\med}{\medskip}
\newcommand{\bp}{\begin{proof}}
\newcommand{\enp}{\end{proof}}

\def\hw{highest weight }

\begin{document}

\title[Matrices of simple spectrum]{Matrices of simple spectrum in irreducible representations of cyclic extension 
of  simple algebraic groups}
\author{A.E. Zalesski }
\address{Department of Physics, Mathematics and Informatics, National Academy of Sciences of Belarus,
66 Prospekt  Nezavisimosti, Minsk, Belarus}
\email{alexandre.zalesski@gmail.com}

\maketitle

Abstract. Let $H$ be a linear algebraic group whose connected
component $G\neq 1$ is simple,  $H/G$ is cyclic and $Z(H)\subseteq Z(G)$. 
We determine the irreducible projective
representations $\phi$ of $H$ such that $\phi(G)$ is \ir and  $\phi(h)$ has simple spectrum for some $h\in H$.  The latter means that all \ir constituents of the group $\phi( \lan h\ran)$ are of \mult 1. (Here $\lan h\ran$ is the subgroup of $H$ generated by $h$.) This extends an earlier known result for $H=G$.

\section{Introduction}

The aim of this paper is to extend the result of  \cite[\S 6]{GS}, \cite{Zsu} on the \ir \reps of simple algebraic groups whose image contains a matrix of simple spectrum. We assume that the \rep field $F$ is algebraically closed;
in this case, the spectrum of  a square matrix $M$ over $F$ is  the set of its \eis and the spectrum of $M$
is called simple if every \ei has \mult 1. In particular, a matrix with simple spectrum is similar to
a diagonal matrix with all diagonal entries are distinct.

 Let $G$ be a linear algebraic group, $\rho$  a rational \rep and  $g\in G$. If $M=\rho(g)$ has simple spectrum then $g$ is a semisimple element and all weights of $\rho$ are of \mult 1. In fact, this is equivalent for $\rho(G)$
to contain a matrix with simple spectrum. For every simple algebraic group the \ir \reps whose all weights are of \mult 1 have been determined in \cite{Zsu}, extending an earlier results by Seitz \cite[\S 6]{GS}.

Recently R. Guralnick raised a question whether this result can be extended to groups $H=\lan G,\si\ran$,
where $G$ is a simple algebraic group and $\si$ is an outer \au of $G$. We prove

\begin{theo}\label{th1} Let G be a simple \ag in defining characteristic $p\geq 0$, let $\si$ be a non-trivial graph \au of G and $H=\lan G,\si\ran$. Let  V be an  H-module such that $V|_G$ is rational  and \irt
Let $\om$ be the \hw   of $V|_G$.
Suppose that some $h\in H$ has simple spectrum on V.   Then one of the \f holds:

$(1)$  There exists $g\in G$ whose spectrum on V is simple;

$(2)$ $G$ is of type $A_2$, $p\neq 2,3$ and   $\om=p^j(\om_1+\om_2)$ for some integer $j\geq 0;$


$(3)$ $G$ is of type $D_4$, $p= 2$ and   $\om=2^j\om_2$ for some integer $j\geq 0.$

\noindent Moreover, in each case $(2),(3)$ there is an element $h\in (H\setminus G)$ with simple spectrum on $V$.
\end{theo}

It is well known that $|\si|\in\{2,3\}$. If $|\si|=p$ and  $h$ is diagonalizable then $h$ has simple spectrum on $V$
\ii so has $g=h^p$ so (1) holds. So we  assume that $(|\si|,p)=1$. In this case $h\in (H\setminus G)$ is diagonalizable \ii so is $h^{|\si|}$.

Theorem \ref{th1} answers the original question by R. Guralnick, and is sufficient to an application he held in mind.
The assumption on $V|_G$ to be rational is not essential due to a theorem of Borel and Tits \cite[Theorem 42]{St}.
The assumption ``$V|_G$ is irreducible" is essential.
In general, the problem of describing the \ir $FH$-modules $V$ containing a simple spectrum matrix but reducible on  $G$ seems to be more complex. We have the \f observations:

\begin{propo}\label{re1} Let $G,H,\si$ be as in Theorem $\rm \ref{th1}$, $k=|\si|$, and let V be an \ir $FH$-module such that $V|_G$
is reducible. Let $h\in (H\setminus G)$. Then the spectrum of h on V is simple \ii $h^k $ has simple spectrum on every
\ir constituent of $V|_G$. In particular, if the spectrum of h on V is simple then the weights of every \ir
constituent  of V are of \mult $1$.
\end{propo}

The converse  is not true in general.

\med

For finite groups of Lie  type $G(q)$  the \ir \reps (over a field of describing characteristic)
whose image contains a matrix with simple spectrum are determined in \cite{sz98,sz00}.
The problem of extending the results of \cite{sz98,sz00} to cyclic extensions $H(q)=\lan  G(q),h\ran $, where $h\in H$ normalizes $G(q)$,  remains open.  We obtain some partial results in the framework of Theorem \ref{th1}, that is, in the case where $H(q)\subset H$. Then Theorem \ref{th1} reduces
the problem to the cases $(1)-(3)$ of it, and the case (1) is sorted out in \cite{sz98,sz00}.

\begin{theo}\label{qq2} Let $G(q)=SL_3(q)$ or $SU_3(q)$, $(q,6)=1$ and $H(q)=\lan G(q),\si\ran$, where $\si$ is the transpose-inverse \au of $G(q)$. Let $V$ be an \ir $FH(q)$-module such that $V|_{G(q)}$ is \ir with highest weight $\om_1+\om_2$. Then $h$ has simple spectrum on V for some $h\in( H(q)\setminus G(q))$.
\end{theo}

 The sentence ``$V|_{G(q)}$ is \ir with highest weight $\om_1+\om_2$"  means $V|_{G(q)}=\overline{V}|_{G(q)}$,
where $\overline{V}$ is an \ir $SL_3(F)$-module with \hw $\om_1+\om_2$.

Similar results are obtained for $H(q)=\lan D_4(q)\cdot \si\ran $ and $H(q)=\lan {}^3D_4(q)\cdot \si\ran $, where $q>8$ even and $\si$ is a triality \au of $D_4(q)$ or ${}^3D_4(q)$, see Sections 6.1 and 6.2. The cases of $q=4,8$ are still open. Note that the answer is opposite for $q=2$ in both the cases, as for $h\in H(2)$ of odd order we have $|h|<d$, where $d$ is the degree of any faithful $\si$-stable \irr $\phi$ of $G(q)$.

\bigskip
{\it Notation} We keep $G$ to denote a simply connected simple algebraic group defined over an algebraically closed field $F$. All representations and modules are assumed to be rational. For general facts of the theory of algebraic groups and their \reps we refer to \cite{St} and \cite{Bo,Bo8}; in particular, the ordering of 
simple roots of $G$ is  as in \cite{Bo}. If $G$ is of rank $n$ then $\om_1\ld \om_n$ are fundamental weights and $\al_1\ld \al_n$ are simple roots of $G$. If  $V$ is a $G$-module and $X$ is a subgroup of $G$ then $V|_X$ denote the restriction of $V$ to $X$. We write  $|X|$ for the order of $X$, and also $|x|$ for the order of an element $x\in G$; $\lan x \ran$ means the group generated by $x$. The normalizer and the centralizer of $X$ in $G$ are denoted by $N_G(X)$ and $C_G(X)$, respectively. $Z(X)$ denotes the center of a group $X$. By $\ZZ$, $\QQ, $ and $\FF_q$ we denote the ring of integers,  the field of rational numbers and the finite field of $q$ elements, respectively. A diagonal matrix with diagonal entries $d_1\ld d_n$ is denoted by $\diag(d_1\ld d_n)$.

\newpage
\section{Reduction to \reps with non-zero weights  of \\ \mult $1$}

The following lemma is obvious:

\bl{ty2} Let $h\in GL_n(F)$ be a diagonalizable element and $l\geq 1$ is an integer. If h has simple spectrum then the \mult of every \ei of $h^l$ does not exceed $l$.\el

\bl{wk1}  Let $M$ be a vector space over an \acf and let $M=M_1\oplus...\oplus M_l$, where $d:=\dim M_1>1$. Let $x\in GL(M)$ be a semisimple element; suppose that $xM_l=M_1$, $xM_i=M_{i+1}$ for $i=1\ld l-1$, and that $x^l$
acts scalarly on $M_1$. Then each \ei of $x$ on $M$ is of \mult $d$.\el

\bp This is straightforward.\enp

In the following statement we do not assume $V|_G$ to be \irt

\begin{propo}\label{gu1} Let $G,H,\si$ be as in Theorem {\rm \ref{th1}}. 
Let V be an H-module such that $V|_G$ is rational and $h\in H$. Suppose that 
h has simple spectrum on V. Then all non-zero weights  of G on V are of \mult $1$ and the \mult of weight $0$ is at most $k:=|\si|$.
\end{propo}

\bp If $h\in G$ then $h\in T$ for some maximal torus of $G$. Then obviously, all $T$-weights are of multiplicity 1, as required. So we assume $h\notin G$. Set  $g=h^k$ so we have $g\in G$.
Then $g\in G$ is semisimple, and hence $g\in T$, where $T$ is some maximal torus of $G$. Moreover,  $T$ can be chosen $h$-stable. (Indeed, $C_G(g)^0$, the connected component of $C_G(g)$, is an $h$-stable semisimple algebraic group on which $h$ acts as a semisimple automorphism; by \cite[Theorem 7.5]{st2},
$C_G(g)^0$ has an $h$-stable maximal torus, which remains maximal in $G$.) We fix this $T$, so $h\in N_H(T)$.)

  Set $W_H=N_H(T)/T$
and $W=N_G(T)/T$. Then $W$ is the Weyl group of $G$ and $W_H=\lan \si,W\ran$. The conjugation action of $N_H(T)$ on $T$ induces the action of $W_H$ on $\Omega:={\rm Hom}\,(T,F^\times)\cong \ZZ^n$, the set of $T$-weights of $G$. Denote by $\om(V)$  the set of $T$-weights of $V$. Then $\om(V)$ is stable under $W_H$. For $\lam\in\om(V)$ denote by $V_\lam$ the $\lam$-weight space of $V$. Then  $hV_\lam=V_{\si(\lam)}$ for every $\lam\in\Omega$, so $h$ permutes $V_\lam$'s. (Here we define $\si(\lam)(t)=\lam(\si(t))$ and $\si(t)=h\up th$ for $t\in T$.)
Clearly,  $g$ stabilizes every $V_\lam$ and acts scalarly on it.

Suppose the contrary, and choose $0\neq \lam\in \om(V)$ so that $\dim V_\lam>1$. Then $y(\lam)\in \om(V)$
and $\dim V_{y(\lam)}=\dim V_\lam$ for every $y\in W_H$.

\medskip

(i) Suppose first that  $hV_{w(\lam)}\neq V_{w(\lam)}$ for some $w\in W$. 
  Since $k$ is a prime, $\sum_{i=0}^{k-1} h^iV_{w(\lam)}$ is a direct sum as
  any sum of distinct weight spaces is a direct sum of them.
By Lemma \ref{wk1}, $h$ has an \ei of \mult at least 2 (in fact,
the \mult in question equals $\dim V_\lam$).
(Indeed, set $M_1=V_{w(\lam)}$, $M_i=h^iV_{w(\lam)}$ in Lemma \ref{wk1}.) So we have a contradiction.

\medskip
(ii) Suppose that $hV_{w(\lam)}= V_{w(\lam)}$ for all $w\in W$.
Then $V_{w(\lam)}=hV_{w(\lam)}= V_{\si(w(\lam))}$ and hence
$\si$ fixes $w(\lam)$ or all $w\in W$.


Next let $R$ be the vector space over $\QQ$ spanned by $\Omega.$ Then $W_H\subset GL(R)$, and $W$ is well known to be a (finite) \ir subgroup of $ GL(R)$. Let $L$ be the $\QQ$-span of the set $\{w(\lam):w\in W\}$. Then $L=R$ as $WL=L$,  $L\neq 0$ and  $W$ is \irt This means that $\si$ fixes every weight of $G$, and hence $\si(t)=t$ for every $t\in T$. In addition, $\si(\al)=\al$ for every root $\al$ of $G$ (as $\al\in R$).
This implies that $\si$ is a trivial \au of $G$, which is a contradiction.


Let $V_0=\{v\in V: tv=v$ for all $t\in T\}$, so  $V_0$ is the zero weight space on $V$. Then 
  $hV_0=V_0$ and $h^k$ acts on $V_0$ trivially. Let $h_0$ be the restriction of $h$ to $V_0$. Then $h_0$ has simple spectrum on $V_0$, and the claim follows from Lemma \ref{ty2}.\enp

\med
Remark. In Proposition \ref{gu1},  If $V|_G$ is irreducible and tensor-decomposable then the weight 0 has \mult 0 or 1. This follows from \cite[Theorem 2(2)]{TZ2}.

\section{Proof of Theorem \ref{th1}}

Below $F$ is an \acf of characteristic $p\geq 0$ and $F^\times$ is the multiplicative group of $F$. Note that
the simple algebraic groups of types $A_1,B_n,n>1,C_n,n>1,E_7,E_8,F_4,G_2$  have no non-trivial graph automorphism.
It is well known that $|\si|=2$ for groups $G$ of type $A_n$, $n>1$, $D_n,$ $n>4 $ and $E_6$, and $|\si|\in\{2,3\}$ for $G$ of type $D_4$. 

\bl{ss1}  Let $G=SL_n(F)$, $p\neq 2$, let  $\si$ be the transpose-inverse \au of $G$ and $H=\lan G,\si\ran$. Let T be the group of diagonal matrices in G and let $h\in (H\setminus G)$. Suppose that $h^2\neq 1$ is a semisimple element of G.

 $(1)$ There is a conjugate $h_1$ of h such that $h_1(T)=T$ and $h_1^2\in T$.

 $(2)$ Let $\overline {h}_1$, $\overline {\si}$ be the projections of $h_1,\si$ into $W_H=N_H(T)/T$. Then $\overline {h}_1=w\overline {\si}$, where $1\neq w\in W=N_G(T)/T$,  $w^2=1$ and $w\overline {\si}=\overline {\si} w$.
\el

\bp We have seen in the proof of Lemma \ref{gu1} that there exists a maximal torus $T_1$ of $G$ such that $h(T_1)=T_1$ and $h^2\in T_1$. As the
maximal tori of $G$ are conjugate, we have   $gT_1g\up=T$ for some $g\in G$, and then $h_1=ghg\up$ is as required in (1).

It is well known that $H=\lan G,\si\ran$. As $\si(t)=t\up$ for every $t\in T$, it follows that $\si$ acts as $-\Id$ on $\Omega={\rm Hom}(T,F^\times)$.
So $\si w=w\si$ for every $w\in W$. 
As $W_H=\lan W,\overline {\si}\ran$, we have $\overline {h}_1=w\overline {\si}$ for some $w\in W$. Then $\overline {h}_1^2=1$ implies $w^2=1$. If $w=1$
then $h_1=t\si$ for some $t\in T$, whence $h_1^2=t\si t\si=t\cdot t\up=1$ and $h^2=1$, contrary the assumption. So (2) follows.\enp

In view of Proposition \ref{gu1}, to prove Theorem \ref{th1} we have to examine the \ir \rep of $G$ whose all non-zero weights are of \mult 1. As above, let $V_0$ be the zero weight space
in $V$.
 If $\dim V_0\leq 1$ then  all weight multiplicities are equal to 1; in this case there exists $h\in G$ with simple spectrum (see for instance \cite[Proposition 1]{SZ}). So it suffices to deal with the \ir \reps of simple algebraic groups whose all non-zero weights are of \mult 1 and $\dim V_0> 1$. Such \reps  are determined in \cite{TZ2}, in particular,
 these are tensor-indecomposable, and hence  are Frobenius twists of those listed in \cite[Table 2]{TZ2}.  This table is reproduced as Table 1 below; it also provides $\dim V_0$.  By Proposition \ref{gu1}, we ignore the cases with $\dim V_0> |\si| $.  In addition, if $V|_G$ is \ir then
the \hw of it is $\si$-invariant, that is, it remains unchanged under the action of the graph  \au $\tau$, say, of order $k=|\si|$, see \cite[Theorem 1.15.2]{GL}. In terms of highest weights, if $\sum a_i\om_i$ ($a_i\in\ZZ$) is the highest weight of $V|_G$ and $\tau(i)=j$
for the nodes $i,j$ at the Dynkin diagram then $\tau(\sum a_i\om_i)=\sum a_i\om_{\tau(i)}$.

So we have to delete  the entries of Table 1 that do not satisfy the above conditions.

\medskip
Let $G=A_n,n>1$. In this case $|H:G|=2=|\si|$, so we drop the cases with $p=2$ and those with  $\dim V_0>2$.  By Table 1, this implies $n=2$, $\om=\om_1+\om_n$, $p\neq 2,3$ or $n=3$, $\om=2\om_2$, $p\neq 2,3$.

Let $G=D_n,n\geq4$. If $n>4$ then  $|H:G|=2=|\si|$ and there is no entry with $\dim V_0\leq2$. Let $n=4$.
If $|\si|=2 $ then $p\neq 2$, and then $\om=\om_2$ or $2\om_2$ in Table 1. The latter
    contradicts Proposition \ref{gu1}. So  $|a|=3 $. Then $\om=\om_2$ by Table 1, and  Proposition \ref{gu1} implies $p=2$.

Let $G=E_6$. Then    $|H:G|=2=|\si|$ and  $\dim V_0>2$ in Table 1, so we may ignore this case.

Summarizing this we obtain:

\bl{in1} In assumption of Theorem {\rm \ref{th1}} 
one of the \f holds:

$(1)$  There exists $g\in G$ whose spectrum on V is simple;

$(2)$ $G$ is of type $A_2$, $p\neq 2,3$ and   $\om=p^m(\om_1+\om_2);$

$(3)$ $G$ is of type $D_4$, $p= 2$ and   $\om=2^m\om_2;$

$(4)$  $G$ is of type $A_3$, $p\neq 2,3$ and   $\om=2p^m\om_2,$\\
where $m\geq 0$ is an integer.
\el

To prove the main statement of  Theorem \ref{th1} we are left to rule out the case (4).
To prove the additional statement of the theorem, 
we have to construct $h\in H$ with simple spectrum on $V$ in the  cases (2), (3) of Lemma \ref{in1}.

\section{The case $G=A_3$, $\om=2\om_2$, $p\neq 2,3$}

In this section we deal with case (4) of Lemma \ref{in1}.
Let $\om_1,\om_2,\om_3$ be the fundamental weights of $G$.

\bl{ne1} Let $G=SL_4(F)$, $p\neq 2,3$, and let V be an \ir G-module with highest weight $2\om_2$.
Let $H=\lan G,\si\ran$, where $\si$ is the transpose-inverse \au of G. Then V extends to an $FH$-module
and H has no   element with simple spectrum on V.
\el

\bp  Suppose the contrary, and let $h\in H$ have simple spectrum on $V$. Then $h\notin G$ as the zero weight of $V$ has \mult 2 (see Table 1).
Let $T$ be the group of diagonal matrices in $G$.
By Lemma \ref{ss1}, we can assume that $h^2\in T$ and $h=\si n_w t$, where $t\in T$ and $n_w\in N_G(T)$
is such that  $w$, the projection of $n_w$ in $N_G(T)/T$, is an involution. We can choose
$n_w$ up to a multiple in $T$,  and we do this so that  $n_w=\diag(\begin{pmatrix}0&1\\ -1&0\end{pmatrix}, 1,1)$ or
$\diag(\begin{pmatrix}0&1\\ -1&0\end{pmatrix}, \begin{pmatrix}0&1\\ -1&0\end{pmatrix})$.
Then $\si n_w\si=n_w$.

Note that $\si(\mu)=-\mu$ for every weight of $G$, in particular, $\si(2\om_2)=-2\om_2$.
If $0\neq v\in V$ is a vector of weight $2\om_2$ then $t\si(v)=\si t\up v$, so $\si(v)$ lies in the $-2\om_2$-weight
space. In addition, $w(2\om_2)=w(\ep_1+\ep_2-\ep_3-\ep_4)=2\om_2$ \cite[Planche I]{Bo}. So $\si n_w (v)$  lies in the $-2\om_2$-weight space. Therefore, $h$ permutes the  $2\om_2$- and  $-2\om_2$-weight spaces. Let $E$ be the sum of them, so $hE=E$.

Let $t=\diag(t_1,t_2,t_3,(t_1t_2t_3)\up)$. Then
 $h^2=\si n_w t\si n_w t= n_w \si t\si n_w t= n_w t\up n_w t$. This equals to $\diag(-t_1t_2\up,-t_1\up t_2,1,1)$
 or $\diag(-t_1t_2\up,-t_1\up t_2,-t_1\up t_2\up t_3^2,-t_1t_2 t_3^{-2})$, respectively. So $(2\om_2)(h^2)=(\ep_1+\ep_2-\ep_3-\ep_4)(h^2)=1=(-2\om_2)(h^2)$ for both the choice of $w$. So the \eis of $h$
on $E$ are the square roots of 1, that is, 1 and $-1$.

As $h$ stabilizes $V_0=C_V(T)$, the zero weight space of $V$, and $h^2$ acts on $V_0$ trivially, it follows that
the \eis of $h$ on $V_0$ are square roots of unity. So either 1 or $-1$ is an \ei of $h$ of \mult at least 2,
a contradiction.
\enp

\section{The case $G=A_2$, $\om=\om_1+\om_n$, $p\neq 2,3$.}


In this case the module $V=V_\om$ can be realized as the adjoint module for $G=SL_3(F)$, whose underlying space $L$ is the space of $(3\times 3)$-matrices with trace 0, and the action of $G$ on $L$ is the conjugation action $l\ra glg\up$ for $l\in L,g\in G$. Note that $\dim V=8$. We can choose for $T$ the group of   diagonal matrices of $G$ and for $\si$ the transpose-inverse \au of $G$. Then $\si$ acts on $L$
as $l\ra -l^{tr}$, where $l^{tr}$ is the transpose of $l$, see \cite[Ch. VIII, \S 13.1]{Bo8}. For a basis of $L$ we choose the matrix units $E_{ij}$ ($1\leq i\neq j\leq 3)$ and $E_{11}-E_{22},E_{22}-E_{33}$, the two latter span the zero weight space.

Let  $h=\si n_wt$, where $t=\diag(t_1,t_2,(t_1t_2)\up)\in T$ and $n_w=\begin{pmatrix}0&1&0\\ 1&0&0\\0&0&-1 \end{pmatrix}\in N_G(T)$, so $n_w^2=1$ and $\si n_w\si=n_w$. We show that the spectrum of $h $ on $V$ is simple for a suitable choice of $t_1,t_2\in F$. Set $g=h^2$.

Note that $g=\si n_wt\si n_wt=n_w\si t\si n_wt= n_wt\up n_wt=\diag(t_1t_2\up,t_1\up t_2,1)\in T$.
 Note that $E_{ij}$ $(1\leq j,j\leq 3)$ are eigenvectors for $g$, and the \eis of $g$ on $V=L$ are $t_1^2t_2^{-2}$, $t_1^{-2}t_2^2$, each of \mult 1, and 1,$t_1t_2\up$, $t_1\up t_2$, each of \mult 2. If we show that $h$ is not scalar on each of the 2-dimensional $g$-eigenspaces of $V$ then the \eis of $h$ on them are distinct square roots of those of $g$, 
and hence the spectrum of $h$ is simple. 

We have $h(E_{12})=\si n_wt(E_{12})t\up n_w\si =t_1t_2\up \si E_{21}\si=-t_1t_2\up E_{12}$, and similarly, 
 $h(E_{21})=-t_1\up t_2 E_{12}$. For the $g$-eigenspaces of \mult 2 we have
 $h(E_{22}-E_{33})=\si n_wt((E_{22}-E_{33}))t\up n_w\si =-(E_{11}-E_{33}) $, so $h$
is non-scalar on the 1-eigenspace of $h^2$. Similarly,  $h(E_{13})=\si n_wt(E_{13})t\up n_w\si =t_1^2t_2 \si E_{23}\si=-t_1^2t_2 E_{32} $ and $h(E_{23})=\si n_wt((E_{23}))t\up n_w\si =t_1t_2^2 \si E_{13}\si=-t_1t_2^2  E_{31}$. So the \eis of $h$ on $V$ are $-t_1t_2\up, $  $-t_1\up t_2 $, $\pm 1, $  $\pm\sqrt{t_1t_2\up}$,
$\pm\sqrt{t_1\up t_2}$.   Whence the claim.

For $t_2=1$ this implies:

\bl{oo1}
Let $t=\diag(t_1,1, t_1\up )\in T$, $n_w=\begin{pmatrix}0&1&0\\ 1&0&0\\0&0&-1 \end{pmatrix}\in N_G(T)$ and $V,\si$ as above.  Suppose that  $\pm 1$, $-t_1,$ $-t_1\up$, $\pm\sqrt {t_1}$, $\pm\sqrt{t_1\up}$
are distinct. Then $h=\si n_w t$ has simple spectrum on V.
\el

In particular,  if $\zeta^4=-1$, 
$\zeta\in F$ 
and $t_1=\zeta^2$ then  
 these are $\pm 1$, $- \zeta^{\pm2}$, $\pm \zeta^{\pm 1}$, which are all 8-roots of unity. Therefore, for this choice of $h$ the spectrum of $h$ is simple. Note that in the latter case
$t=\diag(\zeta^2,1,\zeta^{-2})$ so $|t|=4$.

\begin{corol}\label{33q} Let $G(q)=SL_3(q)$, $(q,6)=1$ and $H(q)=\lan G(q),\si\ran$, where $\si$ is the transpose-inverse \au of $G(q)$. Let $V$ be an \ir $FH(q)$-module such that $V|_{G(q)}$ is \ir with highest weight $\om_1+\om_2$. Then $h$ has simple spectrum on V for some $h\in (H(q)\setminus G(q))$.
\end{corol}

\bp 
We show that the element $h=\si n_w t$ in Lemma \ref{oo1} can be chosen in $G(q)$. We take $t=\diag(t_1,1,t_1\up)$,
where $t_1\in F$ with $|t_1|=q-1$. As $q\neq 3$ is odd, the elements $1$, $t_1^2$, $t_1$, $t_1\up $ are distinct.
As $t,n_w\in G(q)$, the result follows from Lemma \ref{oo1}.
\enp

Let  $G_1=GL_3(q^2)$. Let $V_1$ be the natural module for $GL_3(q^2)$. We choose a basis $B$ in $V_1$
and realize $G_1$ as a matrix group with respect of $B$.  For $x\in G_1$ we write $\overline{x}$ for
the matrix obtained from $x$ by applying the Galois \au of $\FF_{q^2}/\FF_q$ to every entry of $x$. Keep $\si$
for the transpose-inverse \au of $G_1$. It is well known that $G(q)\cong SU_3(q)$ can be realized as $\{x\in SL_3(q^2):\si(x)=\overline{x}\}$.
Then $\si(G(q))=G(q)$.  Let $T$ be the group
of diagonal matrices in $G(q)$ with respect to this basis. It is easy to check that $T$ is of exponent $q+1$ and  $|T|=(q+1)^2$. For $t\in T$ we can write $t=\diag(t_1,t_2,(t_1t_2)\up)$ with $t_1^{q+1}=t_2^{q+1}=1$.

\begin{corol}\label{uu3} Let $G(q)=SU_3(q)$, $(q,6)=1$ and $H(q)=\lan G(q),\si\ran$, where $\si$ is the transpose-inverse \au of $G(q)$. Let $V$ be an \ir $FH(q)$-module such that $V|_{G(q)}$ is \ir with highest weight $\om_1+\om_2$. Then $h$ has simple spectrum on V for some $h\in (H(q)\setminus G(q))$.
\end{corol}

\bp  We show that the element $h=\si w t$ in Lemma \ref{oo1} can be chosen in $G(q)$. We take $t_2=1$, $|t_1|=q+1$ so $t=\diag(t_1,1,t_1\up)$.
As $p\geq  5$, the elements  $\pm 1$, $-t_1,$ $-t_1\up$, $\pm\sqrt {t_1}$, $\pm\sqrt{t_1\up}$
in Lemma \ref{oo1} are distinct.
We choose $n_w$ as in Lemma \ref{oo1}, that is, $n_w=\begin{pmatrix}0&1&0\\ 1&0&0\\0&0&-1 \end{pmatrix}\in N_G(T)$ under basis $B$.
As $t,n_w\in G(q)$, the result follows from Lemma \ref{oo1}.
\enp

Remark. If  $G(q)$ and $V$ are as in Theorem \ref{qq2} and $q$ is a 3-power then there exists $g\in G  $
that has simple spectrum on $V$. This is a special case of \cite[Theorem 0.4]{sz98}.

\section{The case $G=D_4,$ $p=2$, $\om=\om_2$ and $|\si|=3$}

     The \ir $G$-module $V$ with \hw $\om_2$ is a constituent of the adjoint one,
     so $V$ is the only non-trivial \ir constituent of $L(G)$, where $L(G)$ is the Lie algebra of $G$ viewed as a $G$-module. It is well known that $\dim V_{\om_2}=26$.

      Let $T$ be a maximal torus of $G$. Then $T$ is well known to be parameterized as $(t_1,t_2,t_3,t_4)$ with
      $t_1,t_2,t_3,t_4\in F^\times$, so the weight lattice $\Omega={\rm Hom}\,(T,F^\times)$ is of rank 4. Define $\ep_i\in\Om$  by $\ep_i(t_1,t_2,t_3,t_4)=t_i$ for $i=1,2,3,4$. If $\om=\sum a_i\ep_i\in\Om$  $(a_i\in\ZZ)$, then $\om(t_1,t_2,t_3,t_4)=t_1^{a_1}t_2^{a_2}t_3^{a_3}t_4^{a_4}$. If $\om$ is a linear combination of roots then $a_i$'s are integers; this is the case for the weights of $V_{\om_2}$.  We denote by $R$ the $\QQ$-span of $\Omega$,
      so $\dim R=4$. We assume that $\si(T)=T$, so $\si$ acts on  $R$ as a linear transformation and $\si(\Omega)=\Omega$. Note that $\{\ep_i:i=1,2,3,4\}$ is a basis  of $R$.

      We denote by $\Phi$ the set of roots of $G$, which are $\al_1=\ep_1-\ep_2$, $\al_2=\ep_2-\ep_3$, $\al_4=\ep_3-\ep_4$, $\al_4=\ep_3+\ep_4$,
      $\al_1+\al_2$, $\al_2+\al_3$, $\al_2+\al_4$, $ \al_1+\al_2+\al_3$, $ \al_1+\al_2+\al_4$, $\al_2+\al_2+\al_3$, $\al_1+\al_2+\al_3+\al_4$, $\al_1+2\al_2+\al_3+\al_4$,
    and their negations, see \cite[Planche IV]{Bo}. These are non-zero weights of $L(G)$ and of  $V$.
    One can choose a basis of $V$ to be root elements $X_\al$, $\al\in\Phi$, added by two weight zero elements.
    The outer \au $\si$ of $G$ of order 3 acts on the simple roots $\al_i$, $i=1,2,3,4$, by fixing $\al_2$
    and permuting $\al_1,\al_3,\al_4$, and we can assume that $\si(\al_1)=\al_4$,  $\si(\al_3)=\al_1$, $\si(\al_4)=\al_3$, see \cite[Ch. VIII, \S 13.4]{Bo8}. In addition, $X_\al$ can be chosen so that
    $\si(X_\al)=X_{\si(\al)}$ for every $\al\in \Phi$.

 In particular, we have

 \med
 \noindent $\si(X_{\al_1})=X_{\al_4}$, $\si(X_{\al_3})=X_{\al_1}$, $\si(X_{\al_4})=X_{\al_3}$,
  $\si(X_{\al_1+\al_2})=X_{\al_2+\al_4}$,
  $\si(X_{\al_2+\al_4})=X_{\al_2+\al_3}$,

  \noindent $\si(X_{\al_2+\al_3})=X_{\al_1+\al_2}$, $\si(X_{\al_1+\al_2+\al_3})=X_{\al_1+\al_2+\al_4}$,
 $\si(X_{\al_1+\al_2+\al_4})=X_{\al_2+\al_3+\al_4}$, \\ $\si(X_{\al_2+\al_3+\al_4})=X_{\al_1+\al_2+\al_3}$,
 $\si(X_{\al_2})=X_{\al_2}$,
 $\si(X_{\al_1+\al_2+\al_3+\al_4})=X_{\al_1+\al_2+\al_3+\al_4}$, \\ $\si(X_{\al_1+2\al_2+\al_3+\al_4})=X_{\al_1+2\al_2+\al_3+\al_4}$.

\med
The expressions $\al_i=\ep_i-\ep_{i+1}$ for $i=1,2,3$ and $\al_4=\ep_3+\ep_4$  yield the following formulae for the \eis of $t=(t_1,t_2,t_3,t_4)$ on $X_{\al_i}$ for $i=1,2,3,4$:

\begin{equation}\label{eqe}tX_{\al_1}t\up=t_1t_2\up X_{\al_1}, \, tX_{\al_2}t\up=t_3t_3\up X_{\al_2},\, tX_{\al_3}t\up=t_3t_4\up X_{\al_3},\, tX_{\al_4}t\up=t_3t_4 X_{\al_4}.\end{equation}

\noindent In addition, we have $\al_1+\al_2+\al_3+\al_4=\ep_1+\ep_3$, $\al_1+2\al_2+\al_3+\al_4=\ep_1+\ep_2$, whence the \eis of $t$ on the corresponding root elements $X_{\al_1+\al_2+\al_3+\al_4},X_{\al_1+2\al_2+\al_3+\al_4}$ are $t_1t_3$ and $t_1t_2$.

 Let $h=\si t$ for $t=(t_1,t_2,t_3,t_4)\in F$. Then $X_{\al_2}, X_{\al_1+\al_2+\al_3+\al_4},X_{\al_1+2\al_2+\al_3+\al_4}$ are eigenvectors of $h$, as well as $X_{-\al_2}, X_{-\al_1-\al_2-\al_3-\al_4},X_{-\al_1-2\al_2-\al_3-\al_4}$.
 So we have

\med
 \noindent  $h^3(X_{\al_1})=h^2\si t(X_{\al_1})=t_1t_2\up h^2\si(X_{\al_1})=t_1t_2\up h^2(X_{\al_4})=t_1t_2\up h \si t(X_{\al_4})=\\ =t_1t_2\up t_3t_4 h \si(X_{\al_4})=t_1t_2\up t_3t_4 h (X_{\al_3})=t_1t_2\up t_3t_4 \si t (X_{\al_3})=
  t_1t_2\up t_3t_4t_3t_4\up \si(X_{\al_3})=\\ =t_1t_2\up t_3^2  (X_{\al_1})$. Similarly,

 \med\noindent  $h^3(X_{\al_1+\al_2})=h^2\si t(X_{\al_1+\al_2})=
h^2\si t(X_{\ep_1-\ep_3})=t_1t_3\up h^2 h^2\si (X_{\al_1+\al_2})  =t_1t_3\up h^2(X_{\al_2+\al_4})=t_1t_3\up h \si t(X_{\al_2+\al_4})=t_1t_3\up h \si t(X_{\ep_2+\ep_4})=
  t_1t_3\up  t_2t_4   h \si (X_{\al_2+\al_4})=t_1t_3\up  t_2t_4  \si t(X_{\al_2+\al_3})=t_1t_3\up  t_2t_4  \si t(X_{\ep_2-\ep_4})=t_1t_3\up  t_2t_4t_2t_4\up  \si (X_{\al_2+\al_3})=t_1 t_2^2t_3\up (X_{\al_1+\al_2})$, and

 \med\noindent
  $h^3(X_{\al_1+\al_2+\al_3})=h^2\si t(X_{\al_1+\al_2+\al_3})=h^2\si t(X_{\ep_1-\ep_4})=t_1t_4\up h^2 (X_{\al_1+\al_2+\al_4})=\\ t_1t_4\up h \si t (X_{\ep_1+\ep_4})=t_1t_4\up t_1t_4   h \si  (X_{\al_1+\al_2+\al_4})=
 t_1t_4\up t_1t_4   \si t   (X_{\al_2+\al_3+\al_4})=
t_1^2   \si t   (X_{\ep_2+\ep_3})= t_1^2t_2t_3   \si    (X_{\al_1+\al_2+\al_4})=t_1^2t_2t_3       (X_{\al_1+\al_2+\al_3})$.

\med \itf the \eis of $h$ obtained from this are $(t_1t_2\up t_3^2)^{1/3}$
  $(t_1t_2^2 t_3\up)^{1/3}$, $(t_1^2t_2 t_3)^{1/3}$, $\zeta (t_1t_2\up t_3^2)^{1/3}$,
  $\zeta (t_1t_2^2 t_3\up)^{1/3}$, $\zeta (t_1^2t_2 t_3)^{1/3}$, $\zeta ^2(t_1t_2\up t_3^2)^{1/3}$,
  $\zeta ^2(t_1t_2^2 t_3\up)^{1/3}$ and \\ $\zeta ^2(t_1^2t_2 t_3)^{1/3}$, where $\zeta ^3=1\neq \zeta $.
  In addition, a similar action of $h$ on $X_{-\al}$ with $\al$ a positive root above yields the inverses of these  eigenvalues.

  \med
  The basis of the zero weight space $L_0$ of $L(G)$ can be chosen to be $[X_{\al_i},X_{-\al_i}]$ for $i=1\ld 4$ \cite[Ch. VIII, \S 13.4]{Bo8}. Therefore, $\si$ has 3 distinct \eis on it; as the composition series of $L(G)$
  as a $G$-module has two trivial factors and one isomorphic to $V_0$, it follows that the \eis of $\si$
  on $V_0$ are $\zeta,\zeta\up$. 
  \med

  \noindent Altogether the \eis of $h$ on $V$ are $\zeta , \zeta ^2$, $(t_2t_3\up)^{\pm 1}$, $(t_1t_3)^{\pm 1}$, $(t_1t_2)^{\pm 1}$, 
   $\zeta ^i(t_1t_2\up t_3^2)^{\pm 1/3}$,
  $\zeta ^i(t_1t_2^2 t_3\up)^{\pm 1/3}$, $\zeta ^i(t_1^2t_2 t_3)^{\pm1/3}$ for $i=0,1,2$. One observes that there are
  $t_1,t_2,t_3\in F$ such that all these \eis are distinct.
  So we have

  \bl{dd4} Let $t\in T$ be such that $\ep_i(t)=t_i\in F$ and let $\si$ be such that $\si(\al_2)=\al_2$,
   $\si(\al_1)=\al_4$,  $\si(\al_3)=\al_1$, $\si(\al_4)=\al_3$. Let $h=\si t$. Then the \eis of $h$ on V are
   \begin{equation}\label{eq6} \zeta, \zeta ^2, (t_2t_3\up)^{\pm 1},\,\, (t_1t_3)^{\pm 1},\,\, (t_1t_2)^{\pm 1},\,\,  (t_1t_2\up t_3^2)^{\pm 1/3}\zeta ^i,\,\,  (t_1t_2^2 t_3\up)^{\pm 1/3}\zeta ^i,\,\, (t_1^2t_2 t_3)^{\pm1/3}\zeta ^i\end{equation}

  \noindent for $i=0,1,2$,
  where $\zeta$ is a primitive $3$-root of unity in $F$. These are of \mult $1$ for a suitable choice of $t$.
  \el

\subsection{Group $\lan D_4(q),\si\ran$}

Observe that  $V_{\om_2}|_{G(q)}$ is \irt Moreover, all \ir $FG(q)$-module of dimension 26 are Galois conjugate to each other.

\begin{corol}\label{d43} Let $G(q)=D_4(q)$, q even. Suppose that there exist $t_1,t_2, t_3\in \FF_q^\times$ such that
the elements (of $F$) in {\rm (\ref{eq6})} of Lemma {\rm \ref{dd4}}
are distinct.
Then there is an element in $H(q):=\lan G(q),\si\ran$ with simple spectrum on the $2$-modular \irr of dimension $26$.
\end{corol}

\bp It is well known that the elements  $t\in T$ with $t^{q-1}=1$ lies in $G(q)$ (when this is properly 
included into $G$). So the claim follows from Lemma \ref{dd4}. 
\enp

The \f lemma gives a sufficient condition for the assumption of Corollary \ref{d43} to be satisfied.

\bl{d5e} Set $M_1=\{(t_2t_3\up)^{\pm 1}, (t_1t_3)^{\pm 1}, (t_1t_2)^{\pm 1}\}$ and $M_2=\{(t_1t_2\up t_3^2)^{\pm 1}, (t_1t_2^2 t_3\up)^{\pm 1}$, $(t_1^2t_2 t_3)^{\pm1}\}$. The assumption of Corollary {\rm \ref{d43}} is satisfied if there are elements $t_1,t_2,t_3\in \FF_q$ such that  $|M_1|=|M_2|=6$ and either $3|(q-1)$   
or $3\not|(q-1)$ and $s^3\notin M_2$ for every $s\in M_1$.\el

\bp  Let $x\in \FF_q^\times$ with $|x|=q-1$,  and let $\zeta,s\in F$ with $s^3=x, \zeta^3=1\neq \zeta$.
If $3|(q-1)$ then
  $s\notin \FF_q$ whenever $s^3=x$. As $\zeta^i s$ are distinct for  $i=0,1,2$, it follows that    the elements in Corollary \ref{d43} are distinct \ii and $|M_1|=6$ and $|M_2|=6.$
If $3\not| (q-1)$ then  there exists $s\in \FF_q$ with $s^3=x$. As $\zeta s,\zeta^2 s \notin \FF_q$,
the elements in Corollary \ref{d43} are distinct \ii $|M_1|=6$, $|M_2|=6$ and $s^3\notin M_2$ for every $s\in M_1$.\enp

\begin{corol}\label{d48} The assumptions of Corollary {\rm \ref{d43}} can be satisfied for $q>8$.\end{corol}

\bp Let $\xi$ be a primitive $(q-1)$-root of unity. Set $t_1=\xi$,
$t_2=\xi^2$, $t_3=1$. Then $M_1=
\{\xi^{\pm 2},\xi^{\pm 1},\xi^{\pm 3},\}$
and $M_2=\{\xi^{\pm 1},\xi^{\pm 5},\xi^{\pm 4}\}$.
So $|M_1|=6=|M_2|$ for $q>4$. If $3|(q-1)$ then $q>8$ and the result follows by Lemma \ref{d5e}.

Let $(3,q-1)=1$.  Set $M_3=\{s^3:s\in M_1\}$ so $M_3=\{\xi^{\pm 6},\xi^{\pm 3},\xi^{\pm 9} \}$. Then $M_3\cap M_2=\emptyset $ for $q\geq 32$ and the result follows by Lemma \ref{d5e}.\enp

\subsection{Group $\lan {}^3D_4(q)\cdot \si\ran$}

Let $G$ be a simple \ag of type $D_4$ in defining characteristic $p$ and ${\rm Fr}$ a Frobenius endomorphism of $G$. 
This can be chosen so that the group $G^{{\rm Fr}}=\{g\in G: {\rm Fr}(g)=g\}$ is  $D_4(q)$ or ${}^3D_4(q)$ \cite[\S 10]{St}. Let $\si$ be a graph \au of  $G$ of order 3 which permutes
the root subgroups $X_{\al_1},X_{\al_3},X_{\al_4}$ and stabilizes $X_{\al_2}$, see \cite[\S 10]{St}.
Then $\si\cdot {\rm Fr}={\rm Fr}\cdot \si$, so $\si(G(q))=G(q)$. Set $H(q)=\lan G(q),\si\ran$.

\begin{propo} Let $p=2$, let V be an \ir G-module with \hw $\om_2$ and $H=\lan G,\si\ran$. Then V extends to H and remains \ir on $H(q)$. If $q>8$ then there exists $h\in H(q)$ with simple spectrum on V.
\end{propo}

Let $D$ be the group of diagonal matrices in $K:=SL_2(F)$ viewed as a maximal torus of $K$,  and let $U_+=\begin{pmatrix}0&1\\ 0&0\end{pmatrix}$ be a root element
of the Lie algebra of $SL_2(F)$. Let $G_{\al_i}\cong SL_2(F)$ be the $A_1$-type subgroup, $i=1,2,3,4$. Then there exists an algebraic group homomorphism
 $\eta_i:SL_2(F)\ra G_{\al_i}$ which yields  Lie algebra homomorphism
 $\eta_i':L(K)\ra  L(G)$ compatible with the action of these groups on their Lie algebras (that is, $\eta_i'$ is a $K$-module \ho $ L(K)\ra  L(G))|_K$). In particular, $\eta_i'(U_+)=X_{\al_i}$. Let  $x=\diag(y,y\up)\in K$. Then $x U_+ x\up =y^2U_+$ and hence $\eta_i(x)(X_{\al_i})\eta_i(x\up)=y^2X_{\al_i}$. Denote by
 $y_i$ the eigenvalue of $x$ in its action on $X_{\al_i}$.
 This allows us to express $y_1,y_2,y_3,y_4$ in terms of $t_1,t_2,t_3,t_4$ in (\ref{eqe}). We have
 \begin{equation}\label{eqe2}y_1^2=t_1t_2\up,\,\, y_2^2=t_2t_3\up,\,\, y_3^2=t_3t_4\up,\,\, y_4^2=t_3t_4,\end{equation}
whence
\begin{equation}\label{eqe3}t_1=y_1^2y_2^2y_3y_4,\,\,t_2=y_2^2y_3\up y_4\up,\,\,t_3=y_3y_4,\,\,t_4=y_3\up y_4  .\end{equation}

    \begin{lemma}\label{d3d} Let $G(q)={}^3D_4(q)$,
    $q\geq 16$ even, and  $H(q)=\lan G(q),\si\ran$, where $\si$ is as above. 
    Let $V$ be an \ir $FH(q)$-module such that $V|_{G(q)}$ is \ir with highest weight $\om_2$. Then $h$ has simple spectrum on V for some $h\in (H(q)\setminus G(q))$.
\end{lemma}

\bp By (\ref{eqe3})
we have
 $$t_3^2=y_3^2y_4^2,\,\,\, t_4^2=y_3^{-2}y_4^2,\,\,\, t_1t_3=y_1^2y_2^2y_3^2y_4^2,\,\,\,  t_2^2=
 y_2^4y_3^2y_4^2,\,\,\, t_1t_2=y_1^2y_2^4y_3^2y_4^2,$$
so

 $$t_1t_2\up t_3^2=y_1^2y_3^2y_4^2,\,\,\,  t_1t_2^2 t_3\up
 =y_1^2y_2^6y_3^2y_4^2,\,\,\, t_1^2t_2t_3=y_1^4y_2^6y_3^4y_4^4 .$$

\noindent
Suppose that $y_3=y_1^q$, $y_4=y_1^{q^2}$, where $|y_1|=q^3-1$. Then we have

$$t_2t_3\up=y_2^2,\,\, t_1t_3=y_1^{2(1+q+q^2)} y_2^2,\,\, t_1t_2=y_1^{2(1+q+q^2)}y_2^4,$$ $$t_1t_2\up t_3^2=y_1^{2(1+q+q^2)},\,\,t_1t_2^2 t_3\up =y_1^{2(1+q+q^2)}y_2^6,\,\, t_1^2t_2t_3=y_1^{4(1+q+q^2)}y_2^6.$$

Set $u=y_1^{2(1+q+q^2)}$. Then $|u|=q-1$, $u\in \FF_q$ as $(q^3-1)/(q-1)=q^2+q+1$. Then we have

$$t_2t_3\up=y_2^2,\,\,\, t_1t_3=u y_2^2,\,\,\, t_1t_2=uy_2^4,\,
t_1t_2\up t_3^2=u,\,\,\, t_1t_2^2 t_3\up =uy_2^6,\,\,\, t_1^2t_2t_3=u^2y_2^6.$$

\noindent Therefore, the \eis of $h=t\si$ on $V$ are

\med
\centerline{$ \zeta ^{\pm 1}\,\,$, $y_2^{\pm 2}\,\,$, $(u y_2^2)^{\pm 1}\,\,$, $(uy_2^4)^{\pm 1}\,\,$,  $u^{\pm 1/3}\zeta ^i\,\,$,  $(uy_2^6 )^{\pm 1/3}\zeta ^i\,\,$, $(u^2y_2^6)^{\pm1/3}\zeta ^i\,\,$ for $i=0,1,2$.}

\med Note that the element $(y_1,y_2,y_1^q,y_1^{q^2})\in T$
lies in $G(q)$ as $ {\rm Fr}(y_1,y_2,y_1^q,y_1^{q^2})=\\
\si\cdot {\rm Fr}_q (y_1,y_2,y_1^q,y_1^{q^2})=
\si(y_1^{q},y_2,y_1^{q^2},y_1^{q^3})=\si(y_1^{q},y_2,y_1^{q^2},y_1)=
(y_1,y_2,y_1^{q},y_1^{q^2})$.
Here ${\rm Fr}_q$ is an endomorphism of $G$ arising from the mapping $x\ra x^q$ for $x\in F$.

\med

If $(3,q-1)=1$ then there exists $y_2\in \FF_q$ such that $u=y_2^3$. Then the \eis are
as follows, where we use $y$ in place of $y_2$:

$$ \zeta ^{\pm 1},\,\, y^{\pm 2},\,\, y^{\pm 5},\,\, y^{\pm 7},\,\, y^{\pm 1}\zeta ^i,\,\, y ^{\pm 3}\zeta ^i,\,\, y^{\pm 4}\zeta ^i\,\,\, {\rm for}\,\,\, i=0,1,2.$$

If $q\geq 32$ and $(3,q-1)=1$ then  $y^{\pm 2}$,  $y^{\pm 5}$, $y^{\pm 7}$,  $y^{\pm 1}$, $y^{\pm 3}$, $y^{\pm 4}$ are distinct, and hence the above \eis are distinct too.

If 
$3|(q-1)$ then we choose $y_2$ with $y_2^2=u$, and the \eis are

$$ \zeta ^{\pm 1},\,\, y^{\pm 2},\,\,  y^{\pm 4},\,\, y^{\pm 6},\,\,  y^{\pm 8/3}\zeta ^i,\,\,  y ^{\pm 10/3}\zeta ^i,\,\, y^{\pm 2/3}\zeta ^i\,\,\, {\rm for}\,\,\, i=0,1,2.$$

If $q\geq 16$ and $3|(q-1)$  then $ \zeta ^{\pm 1}$, $y^{\pm 2}$,  $y^{\pm 4}$, $y^{\pm 6}$ are distinct and lie in $\FF_q$. In addition, $y^{\pm 8/3}\zeta ^i$,  $y ^{\pm 10/3}\zeta ^i$, $y^{\pm 2/3}$ are distinct and do not lie in $\FF_q$. \itf  the \eis above are distinct.

It is clear that that  there exist elements $t=(t_1,t_2,t_3,t_4)\in G(q)$ satisfying (\ref{eqe3}).
Then $u=y_1^{2(1+q+q^2)}=(t_1t_2\up)^{1+q+q^2}$ and the conditions $y_3=y_1^q,y_4=y_1^{q^2}$ gives
$t_3t_4\up=(t_1t_2\up)^q$, $t_3t_4=(t_1t_2\up)^{q^2}$.
To satisfy this, one can take
$t_3^2=t_1^{q+q^2}t_2^{-q-q^2}$, $t_4^2=t_1^{-q+q^2}t_2^{q-q^2}$.
So the lemma follows.
\enp







\section{The representations of $H$ reducible on $G$}

\bp[Proof of Proposition {\rm\ref{re1}}]  Let $V_1$ be  an \ir constituent of $V|_G$. 
For a weight $\mu$ of $V_1$ let $V_\mu$ be the $\mu$-weight space of $V_1$.

The "only if" part. Since $k$ is a prime, $V|_G=V_1\oplus \cdots \oplus V_k$, where $V_1\ld V_k$ are \ir $G$-modules, and  $H$ permutes $V_1\ld V_k$ transitively (Clifford's theorem).   As $h^kV_1=V_1$ and $h^k\in G$, we can apply Lemma \ref{wk1} to conclude that $h$ has an \ei of \mult at least $d=\dim V_\mu$. So $d=1$.

The "if" part.  Set $V^+_\mu=\sum h^iV_\mu$. Then the \mult of $h^k$ on $V^+_\mu$ equals $k$. By Lemma \ref{wk1},
the \mult of the \eis of $h$ on $V^+_\mu$ equals 1. Clearly, $V^+_\mu\cap V^+_{\mu'}=0$ whenever
$\mu'\neq \mu$ is a weight of $V_1$. As the \eis of $h^k$ on $V_\mu$ and $V_{\mu'}$ are distinct,    the result follows. \enp

Suppose that all weights of $V_1$ are of \mult 1. A natural question is whether there always  exists $h\in H$ with $h^k$ with simple spectrum on $V_1$. The answer is negative.

Indeed, let $G=SL_4(F)$, where $F$ is an \acf of characteristic not 2. Let $V_1=V_{2\om_1}$ be an \ir $G $-module with \hw $2\om_1$ and $V=V_1\uparrow H$, the induced module. Then $V|_G=V_1+V_2$, where $V_2$ is \ir of \hw $2\om_3$. We show that no element $h\in H$ has simple spectrum on $V$.

As in Section 4 we can express $h$ as $\si n_w t$, where $T$ is the group of diagonal matrices in $G$, $t\in T$ and $n_w=\diag(\begin{pmatrix}0&1\\ -1&0\end{pmatrix}, 1,1)$ or
$\diag(\begin{pmatrix}0&1\\ -1&0\end{pmatrix}, \begin{pmatrix}0&1\\ -1&0\end{pmatrix})$.
Then  $h^2=\diag(x,x\up,1,1)$ or $\diag(x,x\up,y,y\up)$ for some $ x,y\in\F^\times$, see Section 4.
Note that $V_{2\om_1}$ has  weights $\pm \om_2$.
As $\om_2=\frac{1}{2}(\ep_1+\ep_2-\ep_3-\ep_4)$ \cite[Planche 1]{Bo}, we have $\om_2(h^2)=\pm 1$, and hence $\om_2(h^2)=(-\om_2)(h^2)$.  So the spectrum of $h^2$ on $V_{2\om_1}$ is not simple, and hence, by   Lemma \ref{re1},  the spectrum of $h$ on $V$ is not simple.

\med
Acknowledgement. I am indebted to R. Guralnick for rising the problem discussed in this paper and his encouraging comments.

\newpage


\begin{center}
\begin{table}
$$\begin{array}{|l|c|c|}
\hline
~~~~~~~~~~G& {\rm highest ~weight~of}~ V
&\mbox{weight } 0 \mbox{ \mult}\\
\hline
A_n,\, n>1,~p\,\not|~ (n+1)&\om_1+\om_n& n\\
(n,p)\neq (2,3),\, p\,|\,(n+1)&\om_1+\om_n& n-1 \\
A_3,~ p>3&2\om_2& 2     \\
 \hline
B_n, n>2,p\neq 2& \om_2&  n \\
~~~~~p\,\,\,|\,(2n+1)&2\om_1&n  \\
~~~~~p\not|\,(2n+1)&2\om_1& n+1\\
\hline
C_2,\,\,p\neq 2&2\om_1&2 \\
C_2,\,\,p\neq 2,5&2\om_2&2\\
C_3,~p=3&2\om_1&3 \\
C_n,~n>2,\,\,p\neq 2&2\om_1 &n\\
\,\,\,\,\,\,\,\,\,\,(n,p)\neq(3,3),\,\,p\not | n&\om_2& n-1\\
\,\,\,\,\,\,\,\,\,\,(n,p)\neq(3,3),\,\,p\,\,\, | n&\om_2& n-2\\
\hline
D_n,\,\,n>3,\,p>2,\,\, p\,\,\,|\,\,n&2\om_1&  n-2 \\
\,\,\,\,\,\,\,\,\,\,\,\,\,n>3,\,p\neq 2,\,\, p\not|\,n&2\om_1& n-1 \\
 & \om_2& n\\
 \,\,\,\,\,\,\,\,\,\,\,\,\,n>3,\,\,p=2& \om_2&n-(2,n) \\
\hline
E_6\,\,\,\,\,p\neq 3&\om_2&6\\
\,\,\,\,\,\,\,\,\,\,\,\, p=3&\om_2&  5\\
\hline
E_7,\,\,p\neq 2&\om_1 &7 \\
\,\,\,\,\,\,\,\,\,\,\,p=2&\om_1&6\\
 \hline
E_8&\om_8 & 8\\
\hline
F_4,\,\, p\neq 2&\om_1 & 4\\
\,\,\,\,\,\,\,\,\,\,\, p=2&\om_1&  2   \\
\,\,\,\,\,\,\,\,\,\,\, p\neq 3& \om_4 &2 \\
\hline
G_2,\, p\neq 3&\om_2&2\\
\hline
\end{array}$$

\end{table}
\end{center}

\begin{center} Table 1: Irreducible $p$-restricted $G$-modules $V$ whose non-zero
weight spaces are $1$-dimensional and whose zero weight has \mult greater than $1$.
\end{center}

\end{document}